\documentclass{siamltex}
\usepackage{amssymb}
\usepackage{amsmath}
\usepackage{rotating}

\newtheorem{pro}{Proposition}[section]
\newtheorem{them}[pro]{Theorem}

\begin{document}

\title{Generalized Benders Decomposition for one Class of MINLPs with Vector Conic Constraint\footnotemark[1]}
\author{Zhou Wei\footnotemark[2]\ \footnotemark[3] \and M. Montaz Ali\footnotemark[3]\ \footnotemark[4]}
\date{}
\maketitle

\date{}

\begin{abstract} In this paper, we mainly study one class of mixed-integer nonlinear programming problems (MINLPs) with vector conic constraint in Banach spaces. Duality theory of convex vector optimization problems applied to this class of MINLPs is deeply investigated. With the help of duality, we use the generalized Benders decomposition method to establish an algorithm for solving this MINLP. Several convergence theorems on the algorithm are also presented. The convergence theorems generalize and extend the existing results on MINLPs in finite dimension spaces.
\end{abstract}

\renewcommand{\thefootnote}{\fnsymbol{footnote}}

\footnotetext[1]{This research was supported by the National Natural Science Foundations of P. R. China (Grant No. 11401518, No. 11261067 and No. 11371312) and IRTSTYN,  and by the Claude Leon Foundation of South Africa.}

\footnotetext[2]{Department of Mathematics, Yunnan University,
Kunming 650091, People's Republic of China (wzhou@ynu.edu.cn).}

\footnotetext[3]{School of Computer Science and Applied Mathematics, University of the Witwatersrand, Wits 2050, Johannesburg, South Africa.}
%Kunming 650091, P. R. China (wzhou@ynu.edu.cn). }

\footnotetext[4]{TCSE Faculty of Engineering, University of the Witwatersrand, Wits 2050, Johannesburg, South Africa (Montaz.Ali@wits.ac.za).}

\noindent {\bf Key words.} {\it Generalized Benders decomposition; MINLP; Duality theory; Vector optimization}

\noindent
{\bf AMS subject classifications.} {\it 90C11, 90C25, 90C30}

\section{Introduction}

In many optimization problems, decision variables appearing in objective and constraint functions are continuous and discrete. These optimization problems can be modelled as mixed-integer nonlinear programming problems (MINLPs). In general, MINLP is defined mathematically as follows:
\begin{equation}\label{1.1a}
{\rm \mathcal{P}} \left \{
\begin{array}l
\mathop{\rm minimize}\limits_{x,\, y}  \ \ \ f(x, y)\\
 {\rm subject\ to}\ \  g(x, y)\leq 0, \\
\ \ \ \ \ \ \ \ \ \ \ \ \ \ \   x\in X, y\in Y\ {\rm discrete\ variable},
\end{array}
\right.
\end{equation}
where $f: \mathbb{R}^n\times \mathbb{R}^p\rightarrow \mathbb{R}$ and $g: \mathbb{R}^n\times \mathbb{R}^p\rightarrow \mathbb{R}^m$ are nonlinear functions, $X\subset\mathbb{R}^n$, and $Y\subset\mathbb{R}^p$ is a polyhedral set of discrete points.

MINLP problem $\mathcal{P}$ is a natural approach to solve problems by simultaneously optimizing the system structure (discrete) and parameters (continuous). Over the past decades, MINLPs have been used in various applications such as the process industry, chemical engineering design, production planning and control, optimal design of gas or water transmission networks, finance and scheduling problems etc.(cf. \cite{CD,FC,GS1,GS2,TS,TS1} and references therein). Note that two subclasses of mixed-integer linear programming (MILP) problem and nonlinear programming (NLP) problem are embedded in MINLP simultaneously, and thus MINLP problem $\mathcal{P}$ falls into the class of NP-hard problems and becomes one of the most difficult optimization problems. It is known that methods for solving MINLP problem $\mathcal{P}$ mainly fall in two broad classes. One class is heuristic methods which do not provide a guarantee that on termination the incumbent is a minimizer, while the other class is deterministic methods which terminate with a guaranteed solution or an indication that the problem has no integer solution. The deterministic methods for solving MINLP problem $\mathcal{P}$ with convex functions are mainly on NLP/LP based on branch-and-bound method (cf. \cite {L,NV}), extended cutting-plane method (cf. \cite{EM1,WP,WP1}), outer approximation method (cf. \cite{FL,WA1,WA2}), variable and Lagrangean decompositions (cf. \cite{FR,MM}), generalized Benders decomposition (cf. \cite{B,B1}) etc.

Vector optimization relates to functional analysis and mathematical programming, and has been found to play many important roles in economics theory, engineering design, management science, multi-criteria decision making and so on. In recent years, the study on vector optimization has received increasing attentions in the literature (see \cite{B2,F,GJ,J,Lu,M} and references therein). To the best of our knowledge, there is not much literature to study MINLPs in the framework of vector optimization, and from the theoretical viewpoint as well as for applications, it is of significance to continue studying MINLPs in general infinite dimension spaces. Motivated by this, in this paper, we mainly study one class of MINLPs with vector conic constraints in the context of Banach spaces, and aim to establish an appropriate algorithm for solving it. Let $E, Z$ be two Banach spaces, $D$ be a normed linear space and $K$ be a closed convex cone in $Z$ which specifies a partial order $\preceq_K$ on $Z$ as follows:
\begin{equation}\label{1.2a}
  z_1\preceq_K z_2\Longleftrightarrow z_2-z_1\in K\ \ {\rm for\ all\ } z_1, z_2\in Z.
\end{equation}
In this paper, we consider the following MINLP problem (VOP) with vector conic constraint:
\begin{equation}\label{1.3a}
{\rm (VOP)} \left \{
\begin{array}l
\mathop{\rm minimize}\limits_{x,\, y} \ \ \  f(x, y)\\
 {\rm subject\ to}\ \  g(x, y)\preceq_K 0, \\
\ \ \ \ \ \ \ \ \ \ \ \ \ \ \   x\in X, y\in Y\ {\rm discrete\ variable},
\end{array}
\right.
\end{equation}
where $f: E\times D\rightarrow \mathbb{R}$ and $g: E\times D\rightarrow Z$, $X\subset E$ and $Y\subset D$ a set with discrete variables. When we take $E:=\mathbb{R}^n$, $Z:=\mathbb{R}^m$, $D:=\mathbb{R}^p$ and $K:=\mathbb{R}^m_+$, problem (VOP) reduces to the MINLP problem $\mathcal{P}$. Hence it is more general to study this class of MINLPs. With respect to solving problem (VOP), we are inspired by some ideas from generalized Benders decomposition and use this decomposition method to construct an appropriate algorithm for finding the optimal value of problem (VOP).

Benders decomposition was first introduced by Benders \cite{B} and has been applied to a variety of optimization problems such as mixed-integer linear programming, nonlinear programming and MINLPs. It is known that Benders decomposition is an approach for exploiting the structure of mathematical programming problems with complicating variables. Such variables, if temporarily fixed, may render the remaining optimization problem considerably more tractable. For the special class of problems studied by Benders \cite{B}, fixing the complicating variables reduces the given problem to an ordinary linear programming, parameterized by the value of the complicating vectors. Along this line, Geoffrion \cite{Ge1} generalized the Benders decomposition to a broader class of problems where the parameterized subproblem need no longer be a linear programming. Rouhani et al. \cite{RL} and Floudas and Ciric \cite{FC} used the generalized Benders decomposition approach to solve MINLPs which are modelled from practical problems of reactive source planning in power systems and heat exchanger network synthesis respectively.  It is noted that Hooker and Ottosson studied logice-based Benders decomposition, one extension of Benders decomposition, and applied this method to planning and scheduling. Readers are invited to consult references \cite{H1,H2,HO} for more details. For these reasons, the generalized Benders decomposition has been extensively studied by many authors over past decades (cf. \cite{B,B1,Ge,Ge1,Gr,SG} and references therein).

%Recently, this generalized Benders decomposition has been used to solve uncapacitated hub location problems (see \cite{CM,CCL}).

Note that Geoffrion \cite{Ge1} employed the generalized Benders decomposition and nonlinear convex duality (cf. \cite{Ge}) to reformulate MINLP problem $\mathcal{P}$ and derive one equivalent master problem. The algorithm, presented through the generalized Benders decomposition procedure, alternates between solutions of relaxed master problems and nonlinear convex subproblems. Sahinidis and Grossmann \cite{SG} further discussed convergence properties of this generalized Benders decomposition procedure. Inspired by \cite{Ge1,SG}, in this paper, we mainly study the generalized Benders decomposition in vector optimization and use this approach to construct one corresponding algorithm for solving MINLP problem (VOP) of \eqref{1.3a}. To achieve this aim, along the line given by Geoffrion \cite{Ge1}, it is necessary to separate problem (VOP) into many subproblems, establish an equivalent master problem of problem (VOP) and solve the relaxation of master problems. For the equivalence between problem (VOP) and its master problem, we are inspired by Geoffrion \cite{Ge} to study the duality of convex vector optimization problems and proved several duality results (see Section 3).

The paper is organized as follows. In Section 2, we give some definitions and preliminaries used in this paper.  Section 3 is devoted to duality theory results on convex vector optimization problems. Several duality results generalize the corresponding ones obtained in \cite{Ge}. In Section 4, we use the generalized Benders decomposition to establish an algorithm for solving problem (VOP) of \eqref{1.3a} with the help of duality results given in Section 3. The convergence theorems on the algorithm are obtained therein. The conclusion of this paper is presented in Section 5.

\section{Preliminaries} Let $E$ be a Banach space (i.e. complete normed linear space) and $E^*$ denote the dual space of $E$ with dual pairing $\langle \cdot, \cdot\rangle$ between $E^*$ and $E$. Given a set $C\subset E$, let $\overline C$ and ${\rm int}(C)$ denote the {\it norm closure} and the {\it interior} of $C$, respectively. For $x\in E$ and $\delta>0$, denote $B(x, \delta)$ the
{\it open ball} with center $x$ and radius $\delta$.

Let $S$ be a closed convex set of $E$ and $x\in S$. We denote $T(S, x)$ the {\it contingent cone} of $S$ at $x$; that is $v\in T(S, x)$ if and only if there exist a sequence $\{v_k\}$ in $E$ converging to $v$ and a sequence $t_k$ in $(0, +\infty)$ decreasing to $0$ such that $x+t_kv_k\in S$ for all $k\in\mathbb{N}$, where $\mathbb{N}$ denotes the set of all natural numbers. The {\it normal cone} of $S$ at $x$, denoted by $N(S, x)$, is defined as:
\begin{equation}\label{2.1}
  N(S, x):=\{x^*\in E^*: \langle x^*, v\rangle\leq 0\ \ {\rm for\ all} \ v\in T(S, x)\}.
\end{equation}
It is known that $N(S, x)$ and $T(S, x)$ are the dual cones of each other  and one can verify that
\begin{equation}
N(S, x)=\{x^*\in E^*: \langle x^*, y-x\rangle\leq 0\ \ {\rm for\ all} \ y\in S\}.
\end{equation}

Let $\psi: E\rightarrow \mathbb{R}\cup\{+\infty\}$ be a convex function. We denote
$$
{\rm dom}(\psi):=\{x\in E: \psi(x)\in\mathbb{R}\}\,\,\,{\rm and}\,\,\,{\rm epi}(\psi)=\{(x, \alpha)\in E\times \mathbb{R}: \psi(x)\leq \alpha\}
$$
the {\it domain} and the {\it epigraph} of $\psi$, respectively. Recall that $\psi$ is said to be {\it lower semicontinuous} at $x\in E$, if $\liminf_{y\rightarrow x}\psi(y)\geq \psi(x)$.  Let $x\in {\rm dom}(\psi)$.
Recall that the {\it subdifferential}  of $\psi$ at $x$ is defined by
\begin{equation}
\partial\psi(x):=\{x^*\in E^*: (x^*, -1)\in N({\rm epi}(\psi), (x, \psi(x)))\}.
\end{equation}
It is known that $x^*\in\partial\psi(x)$ if and only if
$$
\langle x^*, y-x\rangle\leq\psi(y)-\psi(x)\ \ {\rm for\ all} \ y\in E.
$$

Let $F: E\rightrightarrows E^*$ be a set-valued mapping. We denote $
{\rm dom}(F):=\{x\in E: F(x)\not=\emptyset\}$ the {\it domain} of $F$. Let $x\in {\rm dom}(F)$. Recall that $F$ is said to be {\it norm-to-weak$^*$ upper semicontinuous} at $x$, if for every weak$^*$ open set $V$ containing $F(x)$ and every sequence $\{x_n\}\subset {\rm dom}(F)$ with $\|x_n-x\|\rightarrow 0$, one has $F(x_n)\subset V$ for all sufficiently large $n$. Equivalently, it is easy to verify that $F$ is norm-to-weak$^*$ upper semicontinuous at $x$ if and only if for any generalized sequences $\{x_k\}$ and $\{x_k^*\}$ satisfying $x_k\stackrel{\|\cdot\|}\longrightarrow x$, $x_k^*\stackrel{w^*}\longrightarrow x^*$ and $x_k^*\in F(x_k)$ for all $k$, one has $x^*\in F(x)$. Recall that $F$ is said to be {\it locally bounded} at $x$, if there exist constants $\delta, M\in (0, +\infty)$ such that $\|u^*\|\leq M$ holds for any $u\in B(x, \delta)$ and any $u^*\in F(u)$.

Let $Z$ be a Banach space and $K\subset Z$ be a closed convex cone. The {\it partial order} in $Z$ by $K$ is defined as \eqref{1.2a}.
Let $\varphi: E\rightarrow Z$ be a function. Recall that $\varphi$ is said to be {\it K-convex}, if
$$
\varphi(\lambda x_1+(1-\lambda)x_2)\preceq_K\lambda \varphi(x_1)+(1-\lambda)\varphi(x_2)\ \ {\rm for \ any} \ x_1, x_2\in E \ {\rm and\ any}  \  \lambda\in [0, 1].
$$
If one takes $Z:=\mathbb{R}$ and $K:=[0, +\infty)$, then $K$-convexity of $\varphi$ reduces to the general convexity of real-valued function $\varphi$.

\section{Duality for convex vector optimization problems}
In this section, we study one nonlinear convex primal vector optimization problem as well as its associate duality problem and pay main attention to duality theory for this problem. Note that Geoffrion \cite{Ge} investigated duality theory for nonlinear convex programming with convex function constraints in finite dimension spaces, and gave optimality and weak and strong duality theorems by virtue of the concept of perturbation function. Along the line in \cite{Ge}, we apply this approach to the study on convex vector optimization primal problem in Banach spaces and endeavour to provide several duality results on this primal problem and its dual problem. These duality results on primal problem and its dual will play a key role in the construction of generalized Benders decomposition algorithm for solving problme (VOP) of \eqref{1.2a}. Furthermore the finite convergence of this algorithm is also mainly dependent on these duality results (see Section 4). We begin with this convex vector optimization primal problem.

Let $E, Z$ be two Banach spaces and $K\subset Z$ be a closed convex cone with a nonempty interior. We define the partial order in $Z$ by $K$ as follows: for any $z_1, z_2\in Z$,
$$
z_1\preceq_K z_2\Longleftrightarrow z_2-z_1\in K\ {\rm and} \ z_1\prec_K z_2\Longleftrightarrow z_2-z_1\in{\rm int}(K).
$$
We consider the following convex {\it primal} programming problem:
\begin{equation}\label{3.1a}
{\rm (P)} \left \{
\begin{array}l
\mathop{\rm minimize}\limits_{x} \ \ \ f(x)\\
 {\rm subject\ to}\ \  g(x)\preceq_K 0, \\
\ \ \ \ \ \ \ \ \ \ \ \ \ \ \   x\in X,
\end{array}
\right.
\end{equation}
where $X\subset E$ is convex, $f: X\rightarrow \mathbb{R}$ is convex and $g: X\rightarrow Z$ is $K$-convex.

The {\it dual} problem of (P) is taken to be:
\begin{equation}\label{3.2a}
{\rm (D)} \left \{
\begin{array}l
\mathop{\rm maximize}\limits_{u^*}\ \ \big[\inf\limits_{x\in X} \{f(x)+ \langle u^*, g(x)\rangle\}\big]\\
 {\rm subject\ to}\ \   u^*\in K^+,
\end{array}
\right.
\end{equation}
where $K^+:=\{z^*\in Z^*: \langle z^*, z\rangle\geq 0,\ \forall z\in K\}$ denotes the {\it dual cone} of $K$. If we take $E: =\mathbb{R}^n$, $Z:=\mathbb{R}^m$ and $K:=\mathbb{R}^m_+$, primal problem (P) and dual problem (D) reduce to the classic convex programming with convex function constraints and its associate dual programming in finite dimensional spaces, respectively.

Problems (P) and (D) are in close connection with each other and always have optimal values (possibly $\pm\infty$) provided we invoke the customary convention that an infimum (resp. supremum) taken over an empty set is $+\infty$ (resp. $-\infty$). To investigate the interrelationship between problems (P) and (D), we first recall some definitions on problems (P) and (D).

\begin{definition}
A linear continuous functional $u^*\in K^+$ is said to be {\it essentially infeasible} in problem {\rm(D)}, if
$$
\inf\limits_{x\in X} \{f(x)+ \langle u^*, g(x)\rangle\}=-\infty.
$$
Problem (D) is said to be {\it essentially infeasible}, if every $u^*\in K^+$ is essentially infeasible in problem {\rm(D)}; otherwise, problem (D) is said to be {\it essentially feasible}.
\end{definition}

\begin{definition} A pair $(\bar x, \bar u^*)\in E\times Z^*$ is said to satisfy the {\it optimality conditions} for problem (P), if
\begin{equation}\label{3.3b}
\left\{\begin{array}l
{\rm (i)}\ f(\bar x)+\langle \bar u^*, g(\bar x)\rangle=\min\limits_{x\in X}\big\{f(x)+\langle \bar u^*, g(x)\rangle\big\},\\
{\rm (ii)}\ \langle \bar u^*, g(\bar x)\rangle=0,\\
 {\rm (iii)}\ \bar u^*\in K^+, \\
 {\rm (iv)}\ g(\bar x)\preceq_K0.
\end{array}
\right.
\end{equation}
A linear continuous functional $\bar u^*\in Z^*$ is said to be an {\it optimal Lagrange multiplier} for problem (P), if $(\bar x, \bar u^*)$ satisfies the  optimality conditions for some $\bar x\in X$.
\end{definition}

{\bf Remark 3.1.} (a) It is easy to verify that if $\bar u^*\in Z^*$ is an optimal Lagrange multiplier, then there exists $\bar x\in X$ such that $\bar x$ is an optimal solution to problem (P). This means that an optimal Lagrange multiplier presupposes the existence of an optimal solution to problem (P). Furthermore, if $\bar u^*\in Z^*$ is an optimal Lagrange multiplier, then $(\bar x, \bar u^*)$ satisfies the optimality conditions (i)-(iv) in \eqref{3.3b} for any optimal solution $\bar x$ to problem (P).

(b) The optimality conditions are equivalent to a constrained {\it saddle-point} of the Lagrange function, that is, $(\bar x, \bar u^*)$ satisfies the  optimality conditions (i)-(iv) in \eqref{3.3b} if and only if $(\bar x, \bar u^*)\in X\times K^+$ with $\langle \bar u^*, g(\bar x)\rangle=0$ and
$$
f(\bar x)+\langle u^*, g(\bar x)\rangle\leq f(\bar x)+\langle \bar u^*, g(\bar x)\rangle\leq f(x)+\langle \bar u^*, g(x)\rangle\ \ \forall (x, u^*)\in X\times K^+.
$$

For the case of primal problem (P) when taking $E:=\mathbb{R}^n$, $Z:=\mathbb{R}^m$ and $K:=\mathbb{R}^m_+$, Geoffrion \cite{Ge} exploited the concept of {\it perturbation function} to study convex duality theory between  primal problem (P) and its duality problem (D), and proved the existence of optimal Lagrange multipliers for this primal problem. To delve into the problems (P) and (D) in this section, we consider this notion of perturbation function in vector optimization and study its close interrelationship with optimal Lagrange multipliers for problem (P) of \eqref{3.1a}.

Recall that  the {\it perturbation function $v(\cdot)$} associated with problem (P) is defined on $Z$ as follows:
\begin{equation}\label{3.1}
v: Z\rightarrow \mathbb{R}\cup\{\pm\infty\}, \ z\mapsto v(z):=\inf_{x\in X}\big\{f(x)\ \ {\rm subject \ to }\ g(x)\preceq_K z\big\},
\end{equation}
where each $z\in Z$ is called the {\it perturbation vector} for $v(\cdot)$. We denote
\begin{equation}\label{3.2}
A:=\big\{z\in Z: {\rm there \ exists }\ x\in X\ {\rm such\ that }\ g(x)\preceq_Kz\big\}
\end{equation}
the {\it feasible set} of the perturbed problem. Note that $v(z)=+\infty$ if and only if $z\not\in A$ by the customary convention.

\begin{pro}
{\rm (i)} Let $A$ be defined as \eqref{3.2}. Then $A$ is a convex set and $v(\cdot)$ is convex and monotone nonincreasing on $A$.

{\rm (ii)} Suppose that $X$ is compact, $f$ is lower semicontinuous and that $g$ is continuous. Then $A={\rm dom}(v)$ is a closed subset and $v(\cdot)$ is lower semicontinuous. Furthermore, suppose that the Slater constraint qualification
\begin{equation}\label{3.6a}
g(\hat x)\prec_K0 \ \  {\it for \ some} \ \hat x\in X
\end{equation}
holds. Then $v(\cdot)$ is continuous at $0_Z\in A$.
\end{pro}

{\bf Proof.} (i) The convexity of $A$ and $v(\cdot)$ as well as monotone nonincreasing of $v(\cdot)$ follows from the convexity of $f$ and the $K$-convexity of $g$.

(ii) Since $X$ is compact, $f$ is lower semicontinuous and $g$ is continuous, it is easy to verify that $A$ is closed and for any $\bar z\in A$ there exists $\bar x\in X$ such that $v(\bar z)=f(\bar x)$. Then $A={\rm dom}(v)$ and $v(z)=+\infty$ for all $z\not\in A$. For the lower semicontinuity of $v(\cdot)$, it suffices to prove that $v(\cdot)$ is lower semicontinuous on $A$. Let $z\in A$ and $z_i\rightarrow z$ with $z_i\in A$ for all $i\in\mathbb{N}$. Then for any $i\in\mathbb{N}$, there exists $x_i\in X$ such that $g(x_i)\preceq_K z_i$ and
\begin{equation*}
v(z_i)\geq f(x_i)-\frac{1}{i}.
\end{equation*}
Since $X$ is compact, without loss of generality, we can assume that $x_i\rightarrow x\in X$ (considering subsequence if necessary). It follows that $g(x)\preceq_Kz$ as $K$ is closed and $g$ is continuous. This implies that
$$
\liminf_{i\rightarrow \infty}v(z_i)\geq\liminf_{i\rightarrow \infty}(f(x_i)-\frac{1}{i})\geq f(x)\geq v(z)
$$
as $f$ is lower semicontinuous at $x$.

Now, suppose that Slater constraint qualification \eqref{3.6a} hold. Then $-g(\hat x)\in{\rm int}(K)$ and thus there exists $\delta>0$ such that $-g(\hat x)+B(0_Z, \delta)\subset K$. This implies that $B(0_Z, \delta)\subset A={\rm dom}(v)$. Hence $0_Z\in {\rm int}(A)={\rm int}({\rm dom}(v))$ and it follows from \cite[Proposition 3.3]{P} that $v(\cdot)$ is continuous at $0_Z$. The proof is complete. $\Box$

The following proposition provides an equivalent interpretation to optimal Lagrange multipliers and asserts precisely that the set of optimal Lagrange multipliers is essentially the negative of subdifferential  of perturbation function at the origin.

\begin{pro}
Suppose that problem (P) has an optimal solution and denote $U$ the set of all optimal Lagrange multiplier for problem (P). Then $U=-\partial v(0)$.
\end{pro}

{\bf Proof.} The ``$\subset$" part.  Let $\bar u^*\in U$. Then there exists $\bar x\in X$ such that the pair $(\bar x, \bar u^*)$ satisfies the optimal conditions (i)-(iv) in \eqref{3.3b}. From the optimal conditions (i) and (ii), we have
\begin{equation}\label{3.3}
  f(x)+\langle \bar u^*, g(x)\rangle\geq f(\bar x)+\langle \bar u^*, g(\bar x)\rangle=f(\bar x),\ \ \forall x\in X.
\end{equation}
Let $z\in A$ and $x\in X$ with $g(x)\preceq_Kz$. Then $\langle \bar u^*, z\rangle\geq\langle \bar u^*, g(x)\rangle$ by the optimal condition (iii) and it follows from \eqref{3.3} that $f(x)\geq f(\bar x)-\langle \bar u^*, z\rangle$. By taking the infimum of the left-hand side over the indicated values of $x$, one has
\begin{equation*}
  v(z)\geq f(\bar x)-\langle \bar u^*, z\rangle\ \ \forall z\in A.
\end{equation*}
This implies that $-\bar u^*\in \partial v(0)$ since $f(\bar x)=v(0)$ and $v(z)=+\infty$ for any $z\not\in A$.

For the ``$\supset$" part, let $-\bar u^*\in\partial v(0)$ and $\bar x$ be an optimal solution of problem (P). Then $g(\bar x)\preceq_K0$. We only need to prove that the pair $(\bar x, \bar u^*)$ satisfies the optimality conditions (i)-(iii). Noting that $-\bar u^*\in\partial v(0)$, it follows that
\begin{equation}\label{3.4}
v(z)\geq v(0)-\langle\bar u^*, z\rangle\ \ \forall z\in Z.
\end{equation}
This implies that $\langle\bar u^*, z\rangle\geq v(0)-v(z)\geq 0$ holds for all $z\in K$ as $v(\cdot)$ is monotone nonincreasing and consequently $\bar u^*\in K^+$. Noting that decreasing the right-hand side of problem (P) to $g(\bar x)$ will not destroy the optimality of $\bar x$, it follows that
$$
v(g(\bar x))=v(0) \ \ {\rm and} \ \ \langle \bar u^*, g(\bar x)\rangle\geq 0.
$$
On the other hand, $\langle \bar u^*, g(\bar x)\rangle\leq 0$ follows from $-g(\bar x)\in K$ and $\bar u^*\in K^+$. This means that $\langle \bar u^*, g(\bar x)\rangle=0$  and thus the optimality condition (ii) holds. To prove the optimality condition (i), for any $x\in X$, when taking $z:=g(x)$ in \eqref{3.4}, we have
$$
v(g(x))\geq v(0)-\langle \bar u^*, g(x)\rangle.
$$
Since $f(x)\geq v(g(x))$ for all $x\in X$ and $f(\bar x)=v(0)$, one has
$$
f(x)+\langle \bar u^*, g(x)\rangle\geq f(\bar x)=f(\bar x)+\langle \bar u^*, g(\bar x)\rangle\ \ \forall x\in X
$$
(thanks to the optimality condition (ii)). Hence the optimality condition (i) holds. The proof is complete.  $\Box$

It is known from Proposition 3.2 that optimal Lagrange multipliers can be determined from subdifferential $\partial v(0)$ and thus it is necessary to study equivalent conditions ensuring the nonempty of $\partial v(0)$. The following proposition provides a criterion for the existence of subdifferential of perturbation function $v(\cdot)$ at a point where it is finite.

\begin{pro}
Let $v(\cdot)$ associate with problem (P) be defined as \eqref{3.1} and $\bar z\in {\rm dom}(v)$. Suppose that $X$ is compact, $f$ is lower semicontinuous and that $g$ is continuous. Then $\partial v(\bar z)\not=\emptyset$ if and only if there exists $M\in (0, +\infty)$ such that
$$
\frac{v(\bar z)-v(z)}{\|z-\bar z\|}\leq M\ \ \forall z\in Z\backslash\{\bar z\}.
$$
\end{pro}

{\bf Proof.} By virtue of Proposition 3.1, one has $v(\cdot)$ is a lower semicontinuous convex function and $v(z)>-\infty$ for all $z\in Z$. Let $z^*\in \partial v(\bar z)$. Then the necessity part follows by taking $M:=\|z^*\|+1$. It suffices to prove the sufficiency part. Let
\begin{eqnarray*}
\Phi:&=&\{(z, r)\in Z\times \mathbb{R}: v(\bar z)-v(z)\geq r\}\ {\rm and} \\
\Psi: &=&\{(z, r)\in Z\times \mathbb{R}: M\|z-\bar z\|< r\}.
\end{eqnarray*}
Then $\Psi$ and $\Phi$ are convex sets, $\Psi\cap\Phi=\emptyset$ and $\Psi$ is open. By the seperation theorem (cf. \cite[Theorem 3.4]{R}), there exists $(z^*, \beta)\in (Z\times\mathbb{R})^*=Z^*\times\mathbb{R}$ with $(z^*, \beta)\not=(0, 0)$ such that
\begin{equation}\label{3.5}
\sup_{(z, r)\in \Psi}\big\{\langle z^*, z\rangle+\beta r\big\}<\alpha<\inf_{(z, r)\in \Phi}\big\{\langle z^*, z\rangle+\beta r\big\}.
\end{equation}
Then $\beta<0$ (thanks to $(\bar z, 1)\in\Psi$ and $(\bar z, 0)\in\Phi$). Noting that $(\bar z, 0)\in \Phi$ and $(\bar z, \varepsilon)\in\Psi$ for all $\varepsilon>0$, it follows from \eqref{3.5} that
\begin{eqnarray*}
\langle z^*, z-\bar z\rangle+\beta r&<&0\ \ \forall (z, r)\in \Psi \ \ {\rm and} \\
\langle z^*, z-\bar z\rangle+\beta r&\geq& 0\ \ \forall (z, r)\in \Phi.
\end{eqnarray*}
Noting that $(z, v(\bar z)-v(z))\in \Phi$ for all $z\in {\rm dom}(v)$, it follows that
$$
\langle \frac{z^*}{\beta}, z-\bar z\rangle\leq v(z)-v(\bar z),\ \ \forall z\in Z.
$$
This implies that $\widetilde{z}^*:=\frac{z^*}{\beta}\in\partial v(\bar z)$. The proof is complete.  $\Box$

{\bf Remark 3.2.} (a) The proof of Proposition 3.3 is inspired by some ideas from \cite[Lemma 2]{Ge}, and furthermore it is known from the proof that the conclusion is still valid for general proper extended-real-valued, but not taking negative infinity, convex function defined on $Z$.

(b) Under the assumptions of Proposition 3.1(ii), the perturbation function $v(\cdot)$ is a lower semicontinuous convex function. A deep theorem referring to the subdifferential $\partial v$, proved by Br{\o}ndsted and Rockafellar, is that ${\rm dom}(\partial v)$ is dense in ${\rm dom}(v)$; that is for any $z\in Z$ with $v(z)\in\mathbb{R}$, there exists $z_n\rightarrow z$ such that $\partial v(z_n)\not=\emptyset$ for all $n\in \mathbb{N}$. Readers are invited to consult \cite[Theorem 3.17]{P} for Br{\o}ndsted-Rockafellar theorem and its proof in detail.

Using propositions 3.1, 3.2 and 3.3, we obtain the following theorem on characterizations for the existence of optimal Lagrange multipliers.

\begin{them}
Suppose that problem (P) has an optimal solution. Then the following statements are equivalent:

{\rm(i)} The set of all Lagrange multipliers for problem (P) is nonempty.

{\rm(ii)} $\partial v(0)$ is a nonempty set.

{\rm(iii)} $v(0)$ is finite and there exists $M\in (0, +\infty)$ such that
\begin{equation}\label{3.5a}
v(0)-v(z)\leq M\|z\|
\end{equation}
holds for any $z\in Z$.
\end{them}

The following result provides a criterion for the essential feasibility of problem (D) and also gives one necessary condition for essential feasible problem (D).

\begin{pro}
{\rm (i)} Suppose that problem (D) is essentially feasible. Then $v(z)>-\infty$ for all $z\in A$.

{\rm (ii)} Suppose that $X$ is compact, $f$ is lower semicontinuous and that $g$ is continuous. Then problem (D) is essentially feasible.
\end{pro}

{\bf Proof.}
(i) Suppose that problem (D) is essentially feasible. Then there exists $u^*\in K^+$ and $\beta\in \mathbb{R}$ such that
\begin{equation}\label{3.7}
  f(x)+\langle u^*, g(x)\rangle\geq\beta,\ \ \forall x\in X.
\end{equation}
Let $z\in A$ and $x\in X$ with $g(x)\preceq_Kz$. By \eqref{3.7}, one has
$$
f(x)\geq f(x)+\langle u^*, g(x)-z\rangle\geq\beta-\langle u^*, z\rangle.
$$
This implies that $v(z)\geq\beta-\langle u^*, z\rangle>-\infty$.

(ii) By virtue of Proposition 3.1, one has $v(\cdot)$ is a lower semincontinuous convex function and $A={\rm dom}(v)$. Let $z\in A$. Using \cite[Theorem 3.17]{P}, there exists $z_i\rightarrow z$ such that $v(z_i)\in\mathbb{R}$ and $\partial v(z_i)\not=\emptyset$ for all $i$. Then we can take $u_i^*\in\partial v(z_i)$ such that
\begin{equation}\label{3.6}
 v(z)\geq v(z_i)+\langle u_i^*, z-z_i\rangle,\ \ \forall z\in Z,
\end{equation}
and consequently $-u_i^*\in K^+$ due to the nonincreasing of $v(\cdot)$ and \eqref{3.6}. Let $x\in X$. Then $g(x)\in A$ and it follows from \eqref{3.6} that
$$
\langle u_i^*, g(x)-z_i\rangle\leq v(g(x))-v(z_i)\leq f(x)-v(z_i).
$$
This implies that $f(x)+\langle-u_i^*, g(x)\rangle\geq v(z_i)+\langle-u_i^*, z_i\rangle$ and thus
$$
\inf_{x\in X}\big\{f(x)+\langle-u_i^*, g(x)\rangle\}\geq v(z_i)+\langle-u_i^*, z_i\rangle>-\infty.
$$
Hence $-u_i^*$ is feasible to problem (D) and problem (D) is essentially feasible. The proof is complete.  $\Box$

Clearly the customary weak duality result that the optimal value of primal problem (P) is not smaller than the optimal value of problem (D) holds. Furthermore, the next proposition is the strong duality result on problems (P) and (D) which demonstrates the close connection between optimal Lagrange multipliers and solutions to the dual problem (D).
\begin{pro}
Let $v(0)$ be finite. Then $\bar u^*\in Z^*$ is an optimal solution of problem (D) and the optimal values of problem (P) and problem (D) equal if and only if $-\bar u^*\in\partial v(0)$.
\end{pro}

{\bf Proof.} We first prove the sufficiency part. Suppose that $-\bar u^*\in\partial v(0)$. Then
$$
v(z)\geq v(0)+\langle -\bar u^*, z\rangle,\ \ \forall z\in Z.
$$
Thus $\langle \bar u^*, z\rangle\geq v(0)-v(z)\geq 0$ for all $z\in K$ by the nonincreasing of $v(\cdot)$ and $\bar u^*\in K^+$. For any $x\in X$, $f(x)\geq v(g(x))$ and $f(x)+\langle \bar u^*, g(x)\rangle\geq v(0)$. This implies that
\begin{equation}\label{3.13a}
\inf_{x\in X}\{f(x)+\langle \bar u^*, g(x)\rangle\}\geq v(0)=\inf_{x\in X}\{f(x): g(x)\preceq_K 0\}.
\end{equation}
Using the weak duality, one has
\begin{equation}\label{3.14b}
\max_{u^*\in K^+}\Big\{\inf_{x\in X}\{f(x)+\langle u^*, g(x)\rangle\}\Big\}\leq \inf_{x\in X}\{f(x): g(x)\preceq_K 0\}=v(0)
\end{equation}
and it follows from \eqref{3.13a} and \eqref{3.14b} that
$$
\inf_{x\in X}\{f(x)+\langle\bar u^*, g(x)\rangle\}=\max_{u^*\in K^+}\Big\{\inf_{x\in X}\{f(x)+\langle u^*, g(x)\rangle\}\Big\}=v(0).
$$
This means that $\bar u^*$ is an optimal solution of problem (D) and the optimal values of problem (P) and problem (D) equal.

The necessity part. Let $\bar u^*$ be an optimal solution of problem (D) and the optimal values of problems (P) and (D) equal. Then $\bar u^*\in K^+$ and
\begin{equation}\label{3.6b}
\inf_{x\in X}\{f(x)+\langle\bar u^*, g(x)\rangle\}=v(0).
\end{equation}
For any $z\in A$ and any $x\in X$ with $g(x)\preceq_Kz$, by \eqref{3.6b}, one has
$$
f(x)+\langle\bar u^*, z\rangle\geq f(x)+\langle\bar u^*, g(x)\rangle\geq v(0).
$$
This implies that
$$
v(z)\geq v(0)+\langle-\bar u^*, z\rangle \ \ \forall z\in A.
$$
Hence $-\bar u^*\in \partial v(0)$ as $v(z)=+\infty$ for all $z\not\in A$. The proof is complete. $\Box$

We close this section with the following proposition which will be used in next section.
\begin{pro}
Suppose that the optimal value of problem (D) is finite. Then $0_Z\in \overline A$.
\end{pro}

{\bf Proof.} Suppose to the contrary that $0_Z\not\in \overline A$. By the seperation theorem, there exist $z^*\in Z^*$ with $z^*\not=0$ and $\alpha\in \mathbb{R}$ such that
\begin{equation}\label{3.14a}
  \inf_{z\in A}\langle z^*, z\rangle>\alpha>0.
\end{equation}
Let $u^*\in K^+$ such that $u^*$ is essentially feasible in problem (D). Then
$$
\inf_{x\in X}\{f(x)+\langle u^*, g(x)\}>-\infty.
$$
This and \eqref{3.14a} imply that $u^*+tz^*$ is also essentially feasible in problem (D) for all $t>0$ as $g(x)\in A$. Hence
$$
\inf_{x\in X}\{f(x)+\langle u^*+tz^*, g(x)\rangle\}\geq\inf_{x\in X}\{f(x)+\langle u^*, g(x)\rangle\}+t\inf_{x\in X}\langle z^*, g(x)\rangle.
$$
Letting $t\rightarrow +\infty$ and by virtue of \eqref{3.6a}, we obtain the contradiction that optimal value of (D) is $+\infty$. The proof is complete. $\Box$

\section{Generalized Benders decomposition for MINLPs with vector conic constraint}
In this section, we pay main attention to one class of MINLPs with vector conic constraint, and use generalized Benders decomposition approach and duality results on convex vector optimization problems to construct an appropriate algorithm for solving this MINLP. We begin with this class of MINLPs.

Let $E, Z$ be two Banach spaces, $D$ be a normed linear space, and let $X\subset E$ be a closed convex set, $Y\subset D$ be a set with discrete variables and $K\subset Z$ be a closed convex cone with a nonempty interior. This MINLP problem (VOP) is defined as follows
\begin{equation}\label{4.1a}
{\rm (VOP)} \left \{
\begin{array}l
\mathop{\rm  minimize}\limits_{x,\, y} \ \ \ f(x, y)\\
 {\rm subject\ to}\ \  g(x, y)\preceq_K 0, \\
\ \ \ \ \ \ \ \ \ \ \ \ \ \ \   x\in X, y\in Y\ {\rm discrete\ variable},
\end{array}
\right.
\end{equation}
where $f: X\times Y\rightarrow \mathbb{R}$ and $g: X\times Y\rightarrow Z$ satisfy that $f(\cdot, y)$ is convex and $g(\cdot, y)$ is $K$-convex on $X$ for any fixed $y\in Y$.

As pointed out in \cite{Ge1}, Geoffrion employed nonlinear duality theory and generalized Benders decomposition to project MINLP problem  $\mathcal{P}$ in \eqref{1.1a} onto $y$-space, rather than $xy$-space, reformulate problem  $\mathcal{P}$ as one equivalent master problem and construct the generalized Benders decomposition procedure for solving relaxation of master problems. To solve problem (VOP) along this line, we are inspired to separate problem (VOP) into many independent vector optimization problems by fixing discrete variables $y$.

Let $y\in Y$ be fixed. We consider the following vector optimization problem $P(y)$
\begin{equation*}
P(y) \left \{
\begin{array}l
\mathop{\rm minimize}\limits_{x} \ \ \ f(x, y)\\
 {\rm subject\ to}\ \  g(x, y)\preceq_K 0, \\
\ \ \ \ \ \ \ \ \ \ \ \ \ \ \   x\in X,
\end{array}
\right.
\end{equation*}
and its associate dual is defined as follows:
\begin{equation*}
D(y) \left \{
\begin{array}l
\mathop{\rm maximize}\limits_{u^*}\ \ \big[\inf\limits_{x\in X} \{f(x, y)+ \langle u^*, g(x, y)\rangle\}\big]\\
 {\rm subject\ to}\ \   u^*\in K^+.
\end{array}
\right.
\end{equation*}
The perturbation function $v_y(\cdot)$ associated with problem $P(y)$ is defined by
\begin{equation}\label{4.2a}
v_y(z):=\inf_{x\in X}\{f(x, y): g(x, y)\preceq_K z\}, \ \ \forall z\in Z.
\end{equation}

We denote
\begin{equation}\label{4.1}
V:=\{y\in Y: g(x, y)\preceq_K 0\ \ {\rm for\ some}\ x\in X\}
\end{equation}
the feasible set of all values of $y\in Y$ for which vector optimization problem $P(y)$ is feasible. Then problem (VOP) can be equivalently rewritten as
\begin{equation}\label{4.4a}
\left \{
\begin{array}l
\mathop{\rm minimize}\limits_{y} \ \ \ v_y(0)\\
 {\rm subject\ to}\ \  y\in V.
\end{array}
\right.
\end{equation}

After separating problem (VOP) into many problems $P(y)$, it is necessary to establish the master problem which is equivalent to problem (VOP). The key step for this equivalent reformulation is to provide dual equivalent representation for the optimal value $v_y(0)$ of problem $P(y)$ and feasible set $V$. To achieve it, we first give two propositions on the dual equivalent interpretations of $v_y(0)$ and $V$.

The first proposition shows that feasible set $V$ is represented in terms of the intersection of a collection of regions that contain it.

\begin{pro}
Suppose that $X$ is compact and $g(\cdot, y)$ is continuous for any $y\in Y$. Then $\bar y\in V$ if and only if
\begin{equation}\label{4.2}
  \inf_{x\in X}\langle u^*, g(x, \bar y)\rangle\leq 0 \ \ \forall u^*\in K^+\ {\rm with}\ \|u^*\|=1.
\end{equation}
\end{pro}

{\bf Proof.} The necessity part follows from the definition of $V$ in \eqref{4.1}. We only need to prove the sufficiency part. Let
\begin{equation}\label{4.3}
  A(y):=\{z\in Z: g(x, y)\preceq_K z\ {\rm for\ some}\ x\in X\},\ \ \forall y\in Y.
\end{equation}
Since $X$ is compact and $g(\cdot, y)$ is continuous for any $y\in Y$, it follows that $A(y)$ is closed for any $y\in Y$. By \eqref{4.2}, one has
$$
\sup_{u^*\in K^+, \|u^*\|=1}\inf_{x\in X}\langle u^*, g(x, \bar y)\rangle\leq 0.
$$
This implies that
\begin{equation}\label{4.4}
\max_{u^*\in K^+}\Big\{\inf_{x\in X}\langle u^*, g(x, \bar y)\rangle\Big\}=0.
\end{equation}
We consider the following vector optimization problem:
\begin{equation*}
\widetilde{P}(\bar y)\left \{
\begin{array}l
\mathop{\rm minimize}\limits_{x} \ \ \widetilde{f}_{\bar y}(x):=\langle 0_{E^*}, x\rangle\\
 {\rm subject\ to}\ \ g(x, \bar y)\preceq_K 0_Z,\\
\ \ \ \ \ \ \ \ \ \ \ \ \ \ \   x\in X.
\end{array}
\right.
\end{equation*}
Then the dual of problem $\widetilde{P}(\bar y)$ is as follows:
\begin{equation*}
\widetilde{D}(\bar y) \left \{
\begin{array}l
\mathop{\rm maximize}\limits_{u^*}\ \ \big[\inf\limits_{x\in X} \{\langle u^*, g(x, \bar y)\rangle\}\big]\\
 {\rm subject\ to}\ \   u^*\in K^+,
\end{array}
\right.
\end{equation*}
and \eqref{4.4} implies the optimal value of problem $\widetilde{D}(\bar y)$ equals $0$. Using the proof of Proposition 3.7, we obtain that $0_Z\in\overline{A(\bar y)}= A(\bar y)$. Then there exists $x\in X$ such that $g(x, \bar y)\preceq_K0_Z$ and therefore $\bar y\in V$. The proof is complete. $\Box$

By virtue of Theorem 3.4 and Proposition 3.6, we obtain the following proposition on the dual interpretation of $v_y(0)$ which is given by the pointwise supremum of a collection of functions.

\begin{pro}
Suppose that $v_y(0)$ is finite and problem $P(y)$ possesses a Lagrange multiplier for any $y\in V$. Then the optimal value of problem $P(y)$ equals that of its dual problem $D(y)$ for all $y\in V$; that is,
\begin{equation}
  v_y(0)=\max_{u^*\in K^+}\Big\{\inf_{x\in X}\big\{f(x, y)+\langle u^*, g(x, y)\rangle\big\}\Big\}
\end{equation}
holds  for all $y\in V$.
\end{pro}

Under the assumptions that $X$ is compact, $g(\cdot, y)$ is continuous for any $y\in Y$ and problem $P(y)$ possesses an optimal Lagrange multiplier for any $y\in Y$ where problem $P(y)$ is feasible, by using Propositions 4.1 and 4.2, we obtain that problem (VOP) is equivalent to the following master problem:
\begin{equation}\label{4.9a}
\left\{
\begin{array}l
\mathop{\rm minimize}\limits_{y\in Y}\ \ \Big[\sup\limits_{u^*\in K^+}\big\{\inf\limits_{x\in X}\{f(x, y)+\langle u^*, g(x, y)\rangle\}\big\}\Big]\\
{\rm subject \ to}\ \   \inf\limits_{x\in X}\langle z^*, g(x,  y)\rangle\leq 0, \ \ \forall z^*\in K^+\ {\rm with}\ \|z^*\|=1.
\end{array}
\right.
\end{equation}
Using the definition of supremum as the smallest upper bound, the master problem \eqref{4.9a} is equivalent to the following master problem (MP):
\begin{equation}\label{4.7}
{\rm(MP)}\left\{
\begin{array}l
  \mathop{\rm minimize}\limits_{y\in Y,\, \eta\in\mathbb{R}}\ \ \ \eta\\
   {\rm subject\ to}\ \ \eta\geq\inf\limits_{x\in X}\{f(x, y)+\langle u^*, g(x, y)\rangle\}, \ \ \forall u^*\in K^+\\
 \ \ \ \ \ \ \ \ \ \ \ \ \ \ \ \inf\limits_{x\in X}\langle z^*, g(x, y)\rangle\leq 0,\ \ \forall z^*\in K^+\ {\rm with}\ \|z^*\|=1.
\end{array}
\right.
\end{equation}

It is known that one type of relaxation, in which not all constraints are included, is one natural strategy for solving master problem (MP) in \eqref{4.7}. We begin to solve one relaxed version of master problem, not including all constraints in \eqref{4.7}. If the obtained optimal solution does not satisfy constraints having not been considered, then we generate and add to the relaxed problem one or more violated constraints and solve it again. We continue this approach until a relaxed problem solution satisfies all constraints, or until a termination criterion demonstrates that a solution of acceptable accuracy has been obtained. Geoffrion \cite{Ge1} discussed in detail that a solution to a relaxed version of master problem can be tested for feasibility with respect to the ignored constraints and one violated constraint can be generated in case of infeasibility. This discussion given in \cite{Ge1} can also be applied to problem (MP) in \eqref{4.7} and it enables (MP) of \eqref{4.7} to be solved by this relaxation approach.

Now, we can formally state the generalized Benders decomposition procedure for solving problem (VOP). For the validity of equivalence between problems (VOP) and master problem (MP), we suppose that (VOP) satisfies the following assumption:

(A) {\it $X$ is compact, $g(\cdot, y)$ is continuous for any $y\in Y$ and the Slater constraint qualification
 \begin{equation}\label{4.17}
g(\hat x, y)\prec_K0\ \  {\it for \ some} \ \hat x\in X
\end{equation}
holds for any $y\in Y$ where problem $P(y)$ is feasible.}

Using Propositions 3.1 and 3.2, it follows from the Slater constraint qualification \eqref{4.17} that $P(y)$ possesses an optimal Lagrange multiplier for any $y\in Y$ where problem $P(y)$ is feasible.

The detailed algorithm, used to solve problem (VOP) by generalized Benders decomposition procedure, is stated as follows.

{\bf Algorithm 1}(Generalized Benders Decomposition procedure for problem (VOP))

{\bf Step 1}:  Take $y_{1}\in V$ and $z_{1}^*\in K^+$ with $\|z^*_{1}\|=1$. Solve the primal problem
$P(y_{1})$ and obtain an optimal Lagrange multiplier $u^*_{1}$ of $P(y_{1})$. Set
$$T^1=S^1:=\{1\}\ {\rm and} \ {\rm UBD}^1:=v_{y_{1}}(0).
$$
Select the convergence tolerance parameter $\varepsilon>0$ and let $k:=1$.

{\bf Step 2}: Solve the following relaxed master problem ${\rm RMP}(T^k, S^k)$:
\begin{equation}\label{4.8}
{\rm RMP}(T^k, S^k)\left\{
\begin{array}l
  \mathop{\rm minimize}\limits_{y\in Y,\, \eta\in\mathbb{R}} \ \ \ \eta\\  {\rm subject\ to}\ \ \eta\geq\inf\limits_{x\in X}\{f(x, y)+\langle u_{i}^*, g(x, y)\rangle\}, \ \ \forall i\in T^k,\\
 \ \ \ \ \ \ \ \ \ \ \ \ \ \ \  \inf\limits_{x\in X}\langle z_{j}^*, g(x, y)\rangle\leq 0,\ \ \forall j\in S^k.
\end{array}
\right.
\end{equation}
Denote $(y_{k+1}, \eta_{k+1})$ the optimal solution of ${\rm RMP}(T^k, S^k)$. If ${\rm UBD}^k\leq \eta_{k+1}+\varepsilon$, terminate; otherwise, go to {\bf Step 3}.

{\bf Step 3}: Solve the primal problem $P(y_{k+1})$. There must occur one of the following two cases:

{\bf(a)} $v_{y_{k+1}}(0)<+\infty$. If $v_{y_{k+1}}(0)\leq\eta_{k+1}+\varepsilon$, terminate; otherwise, determine an optimal Lagrange multiplier $u_{k+1}^*$ of problem $P(y_{k+1})$, and let
$$
T^{k+1}:=T^k\cup\{k+1\}, S^{k+1}:=S^k\ {\rm and} \ {\rm UBD}^{k+1}:=\min\{{\rm UBD}^{k}, v_{y_{k+1}}(0)\}.
$$
Set $k:=k+1$ and return to {\bf Step 2}.

{\bf(b)} $v_{y_{k+1}}(0)=+\infty$; that is, problem $P(y_{k+1})$ is infeasible. Take $z_{k+1}^*\in K^+$ with $\|z_{k+1}^*\|=1$ such that
\begin{equation*}
  \inf_{x\in X}\langle z_{k+1}^*, g(x, y_{k+1})\rangle>0.
\end{equation*}
Let
$$
T^{k+1}:=T^k, S^{k+1}:=S^k\cup\{k+1\} \ \ {\rm and} \ \ {\rm UBD}^{k+1}:={\rm UBD}^{k}.
$$
Set $k:=k+1$ and return to {\bf Step 2}.

Now, we study the following example and demonstrate the generalized Benders decomposition procedure when solving MINLP problem by Algorithm 1.

{\bf Example 4.1.} Consider the following MINLP problem:
\begin{equation}\label{4.13a}
\left \{
\begin{array}l
\mathop{\rm minimize}\limits_{x,\,y} \ \ \ f(x, y)=-x+\max\{y-1, -y+1\}\\
 {\rm subject\ to}\ \ g(x, y)=x+\max\{-y, y-2\}\leq 0,  \\
\ \ \ \ \ \ \ \ \ \ \ \ \ \ \   x\in [-1, 1],\, y\in \{0, 1, 2\}.
\end{array}
\right.
\end{equation}
Then $X=[-1, 1]$, $Y=\{0,1,2\}$ and $K=K^+=[0,+\infty)$. It is easy to verify that $(x^*, y^*)=(1, 1)$ and $f^*=f(1, 1)$ is the solution of problem \eqref{4.13a}. First, we take $y_1=0$ and $z_1^*=1\in K^+$. Solve primal problem $P(y_1)$ and its dual $D(y_1)$, and denote an optimal Lagrange multiplier $u_1^*=1$. Let $T^1=S^1:=\{1\}, {\rm UBD}^1=v_{y_1}(0)=1$ and $\varepsilon\in (0,1)$. By computing, the relaxed master problem ${\rm RMP}(T^1, S^1)$ is defined as follows:
\begin{equation}
{\rm RMP}(T^1, S^1)\left\{
\begin{array}l
  \mathop{\rm minimize}\limits_{y\in Y,\, \eta\in\mathbb{R}} \ \ \ \eta\\  {\rm subject\ to}\ \ \max\{y-1, -y+1\}+\max\{-y, y-2\}\leq \eta,\\
 \ \ \ \ \ \ \ \ \ \ \ \ \ \ \  -1+\max\{-y, y-2\}\leq 0.
\end{array}
\right.
\end{equation}
It is easy to verify that its solution is $(y_2,\eta_2)=(1,0)$. Since ${\rm UBD}^1>\eta_2+\varepsilon$, by {\bf Step 2}, we go to {\bf Step 3} and solve primal problem $P(y_2)$. Noting that $x_2=1$ is an optimal solution of  $P(y_2)$ and $v_{y_2}(0)=-1<{\rm UBD}^1+\varepsilon$, then terminate the algorithm by {\bf Step 2} and consequently $(x_2,y_2, v_{y_2}(0))=(1,1,-1)$ is an $\varepsilon$-tolerance optimal solution of problem \eqref{4.13a}.

Next, we focus on theoretical convergence of Algorithm 1 by generalized Benders decomposition procedure and prove convergence theorems with the help of some mild assumptions. We first need the following proposition which will be used in the proof of convergence theorems.

\begin{pro}
Suppose that $X$ is compact, and $f(\cdot, \cdot)$, $g(\cdot, \cdot)$ are continuous on $X\times Y$. Denote $U(y)$ the set of all optimal Lagrange multipliers of problem $P(y)$ for any $y\in Y$  where problem $P(y)$ is feasible and let
\begin{equation}\label{4.18}
  L(y, u^*):=\inf_{x\in X}\big\{f(x, y)+\langle u^*, g(x, y)\rangle\big\},\ \ \forall (y, u^*)\in Y\times K^+.
\end{equation}
Then $L(\cdot, \cdot)$ is $\|\cdot\|\times w^*$ continuous on $Y\times K^+$ and the set-valued mapping $U(\cdot)$ is norm-to-weak$^*$ upper semicontinuous on $V$.
\end{pro}

{\bf Proof.} Let $(\bar y, \bar u^*)\in Y\times K^+$, and take any generalized sequence $(y_{\alpha}, u^*_{\alpha})$ in $Y\times K^+$ such that $y_{\alpha}\rightarrow \bar y$ and  $u^*_{\alpha}\stackrel{w^*}\longrightarrow \bar u^*$. Then for any $\alpha$, there exists $x_{\alpha}\in X$ such that
\begin{equation}\label{4.19}
  f(x_{\alpha}, y_{\alpha})+\langle u^*_{\alpha}, g(x_{\alpha}, y_{\alpha})\rangle=L(y_{\alpha}, u^*_{\alpha})
\end{equation}
as $X$ is compact and $f(\cdot, \cdot), g(\cdot, \cdot)$ are continuous. Noting that $X$ is compact, without loss of generality, we can assume that $x_{\alpha}\rightarrow \bar x\in X$ (considering generalized subsequence if necessary). By \eqref{4.19}, one has
\begin{eqnarray*}
L(\bar y, \bar u^*)&\leq& f(\bar x, \bar y)+\langle \bar u^*, g(\bar x, \bar y)\rangle\\&\leq&\lim_{\alpha}(f(x_{\alpha}, y_{\alpha})+\langle u^*_{\alpha}, g(x_{\alpha}, y_{\alpha})\rangle)\\&=&\liminf_{\alpha}L(y_{\alpha}, u^*_{\alpha}).
\end{eqnarray*}
This implies that $L(\cdot, \cdot)$ is $\|\cdot\|\times w^*$ lower semicontinuity at $(\bar y, \bar u^*)$. It suffices to prove the $\|\cdot\|\times w^*$ upper semicontinuity of $L(\cdot, \cdot)$ at $(\bar y, \bar u^*)$. For any $x\in X$, one has
\begin{eqnarray*}
f(x, \bar y)+\langle \bar u^*, g(x, \bar y)\rangle&=&\lim_{\alpha}(f(x, y_{\alpha})+\langle u^*_{\alpha}, g(x, y_{\alpha})\rangle)\\ &\geq&\limsup_{\alpha} L(y_{\alpha}, u_{\alpha}^*),
\end{eqnarray*}
This yields that
$$
L(\bar y, \bar u^*)=\inf_{x\in X}(f(x, \bar y)+\langle \bar u^*, g(x, \bar y)\rangle)\geq\limsup_{\alpha} L(y_{\alpha}, u_{\alpha}^*)
$$
and consequently $L(\cdot, \cdot)$ is $\|\cdot\|\times w^*$ upper semicontinuous at $(\bar y, \bar u^*)$. Thus $L(\cdot, \cdot)$ is $\|\cdot\|\times w^*$ continuous at $(\bar y, \bar u^*)$.

Next, we prove the norm-to-weak$^*$ upper semicontinuity of mapping $U(\cdot)$. Let $\bar y\in V$ and take any generalized sequence $(y_{\alpha}, u^*_{\alpha})$ in $V\times K^+$ such that $y_{\alpha}\rightarrow\bar y$ and $u^*_{\alpha}\stackrel{w^*}\longrightarrow \bar u^*$ with $u^*_{\alpha}\in U(y_{\alpha})$. We only need to show that $\bar u^*\in U(\bar y)$. By Definition 3.2, for any $\alpha$, there exists $x_{\alpha}\in X$ such that $g(x_{\alpha}, y_{\alpha})\preceq_K 0$, $u^*_{\alpha}\in K^+$ with $\langle u^*_{\alpha}, g(x_{\alpha}, y_{\alpha})\rangle=0$ and
\begin{equation}\label{4.20}
  f(x_{\alpha}, y_{\alpha})+\langle u^*_{\alpha}, g(x_{\alpha}, y_{\alpha})\rangle=L(y_{\alpha}, u^*_{\alpha}).
\end{equation}
Since $X$ is compact, without loss of generality, we can assume that $x_{\alpha}\rightarrow \bar x\in X$ (considering generalized subsequence  if necessary) and it follows from the continuity of $g(\cdot, \cdot)$ that $g(\bar x, \bar y)\preceq_K0$, $\langle \bar u^*, g(\bar x, \bar y)\rangle=0$ and $\bar u^*\in K^+$. Using \eqref{4.20} and the continuity of $L(\cdot, \cdot)$ and $f(\cdot, \cdot)$, one has
$$
f(\bar x, \bar y)+\langle \bar u^*, g(\bar x, \bar y)\rangle=L(\bar y, \bar u^*).
$$
This implies that $(\bar u^*, \bar x)$ satisfies the optimality conditions (i)-(iv) in \eqref{3.3b} for $P(\bar y)$ and consequently $\bar u^*\in U(\bar y)$.
The proof is complete. $\Box$

\begin{them}
Suppose that $X, Y$ are compact and $f(\cdot, \cdot)$ and $g(\cdot, \cdot)$ are continuous on $X\times Y$. Denote $U(y)$ the set of all optimal Lagrange multipliers of $P(y)$ for any $y\in Y$ for which problem $P(y)$ is feasible and suppose that the set-valued mapping $y\mapsto U(y)$ is locally bounded on $V$. Then for any given $\varepsilon>0$, the algorithm by generalized Benders decomposition procedure terminates in a finite number of steps.
\end{them}

{\bf Proof.} Suppose to the contrary that there exists $\varepsilon_0>0$ such that the procedure does not terminate in a finite number of steps. Then there exists a sequence $(y_k, \eta_k)$ in $V\times \mathbb{R}$ generated by {\bf Step 2}.  For any $k$, we take $u^*_k\in U(y_k)$. From  {\bf Step 2}, it is not hard to verify that $\{\eta_k\}$ is nondecreasing and bounded above. By taking a generalized subsequence if necessary, we can assume that $(y_k, \eta_k)\rightarrow (\bar y, \bar \eta)\in V\times \mathbb{R}$ since $X, Y$ are compact and $V$ is closed. Noting that mapping $U(\cdot)$ is locally bounded at $\bar y$, it follows that $\{u^*_k\}$ is bounded. Applying Banach-Alaoglu theorem (cf. \cite[Theorem 3.15]{R}), we can assume that $u_k^*\stackrel{w^*}\longrightarrow \bar u^*$ (considering the generalized subsequence if necessary). By the norm-to-weak$^*$ upper semicontinuity of $U(\cdot)$ in Proposition 4.3, one has $\bar u^*\in U(\bar y)$. Using {\bf Step 2} and {\bf Step 3}, we yield
$$
\eta_{k+1}\geq L(y_{k+1}, u_k^*)
$$
and consequently
\begin{equation}\label{4.21}
\bar\eta\geq L(\bar y, \bar u^*)
\end{equation}
by taking limits as $k$. Noting that $u^*_k\in U(y_k)$ and $\bar u^*\in U(\bar y)$, it follows from Propositions 3.2 and 3.6 that
$$
v_{y_k}(0)=L(y_k, u_k^*) \ \  {\rm and} \ \ v_{\bar y}(0)=L(\bar y, \bar u^*).
$$
This and the $\|\cdot\|\times w^*$ continuity of $L(\cdot, \cdot)$ imply that
\begin{equation}\label{4.22}
  \lim_{k}v_{y_k}(0)=v_{\bar y}(0).
\end{equation}
Then, for $\varepsilon_0>0$, when $k$ is sufficiently large, one has
$$
\eta_k+\frac{\varepsilon_0}{2}>\bar\eta\geq v_{\bar y}(0)>v_{y_k}(0)-\frac{\varepsilon_0}{2}
$$
(thanks to \eqref{4.21} and \eqref{4.22}). Thus
$$
v_{y_k}(0)\leq\eta_k+\varepsilon_0\ \ {\rm for\  all}\ k \ {\rm sufficiently\ large},
$$
which contradicts the termination criterion at {\bf Step 3(a)}. The proof is complete. $\Box$

The following convergence theorem can be obtained from Theorem 4.1.

\begin{them}
Suppose that $X$ is compact, the cardinality of $Y$ is finite and that $f(\cdot, y)$ is continuous on $X$ for any fixed $y\in Y$. Then for any given $\varepsilon>0$, the algorithm by generalized Benders decomposition procedure terminates in a finite number of steps.
\end{them}

{\bf Proof.} Let $U(y)$ denote the set of all optimal Lagrange multipliers of problem $P(y)$  for which problem $P(y)$ is feasible. Let $y\in Y$ such that problem $P(y)$ is feasible. By virtue of the Slater constraint qualification \eqref{4.17} and Proposition 3.1, one has $v_y(\cdot)$ is continuous at $0\in Z$ and it follows from Theorem 3.4 that
\begin{equation}\label{4.18a}
  U(y)=-\partial v_y(0).
\end{equation}
Since $v_y(\cdot)$ is continuous at $0$, by using \cite[Proposition 1.11]{P}, one has $\partial v_y(0)$ is bounded. Noting that the cardinality of $Y$ is finite, it follows from \eqref{4.18a} that $U(\cdot)$ is bounded on $V$. Hence the termination criterion in Theorem 4.5 follows from Theorem 4.4. The proof is complete. $\Box$

\section{Conclusions}
This paper is devoted to the study on one class of MINLPs with vector conic constraint in the context of Banach spaces. By using convex primal vector optimization programming and its associated duality results obtained in Section 3, the generalized Benders decomposition method has been used to study MINLP problem (VOP) and establish a corresponding algorithm for solving this problem (see Algorithm 1 in Section 4). With regards to the convergence of the algorithm, it is shown by Theorems 4.4 and 4.5 that the termination criterion in a finite number of steps follows with some mild assumptions. The algorithm extends the generalized Benders decomposition in the sense of solving MINLPs from problem $\mathcal{P}$ in finite dimension space to problem (VOP) in more general Banach space.\\

{\bf Acknowledgement.} The authors are indebted to two anonymous referees for
their comments and suggestions which help us to improve our presentation and draw our attention to the work by Hooker on logic-based Benders decomposition.

\end{document}